\documentclass[a4paper,10pt,reqno]{amsart}
\usepackage[utf8]{inputenc}
\usepackage[T1]{fontenc}
\usepackage{amsthm}
\usepackage{amsmath}
\usepackage{amssymb}
\usepackage[inline]{enumitem}
\usepackage{comment}
\PassOptionsToPackage{hyphens}{url}\usepackage{hyperref}
\usepackage{fancyhdr}
\usepackage{mathrsfs}
\usepackage{stmaryrd}
\usepackage[normalem]{ulem}
\usepackage{xcolor}
\usepackage{tikz-cd}
\usetikzlibrary{arrows}

\usepackage{cite}
\setlist[description]{%
  itemsep=0.05cm,               
  font={\normalfont\textsc}, 
 leftmargin=\parindent,
 labelindent=\parindent
}

\newtheorem{theorem}{Theorem}[section]
\newtheorem{lemma}[theorem]{Lemma}
\newtheorem{proposition}[theorem]{Proposition}
\newtheorem{corollary}[theorem]{Corollary}

\theoremstyle{definition}
\newtheorem{definition}[theorem]{Definition}
\newtheorem{example}[theorem]{Example}
\newtheorem{remark}[theorem]{Remark}

\newtheorem*{conjecture*}{Conjecture}
\newtheorem*{problem*}{Problem}
\newtheorem*{question*}{Question}

\definecolor{blue-url}{RGB}{0,0,100}
\definecolor{red-url}{RGB}{100,0,0}
\definecolor{green-url}{RGB}{0,100,0}
\definecolor{light-yellow}{RGB}{255,255,128}
\definecolor{light-blue}{RGB}{193,255,255}
\definecolor{light-red}{RGB}{239,83,80}

\hypersetup{
	pdftitle={Isomorphisms of reduced power monoids of semilattices},
	pdfauthor={},
	pdfmenubar=false,
	pdffitwindow=true,
	pdfstartview=FitH,
	colorlinks=true,
	linkcolor=blue-url,
	citecolor=green-url,
	urlcolor=red-url
}

\renewcommand{\setminus}{\smallsetminus}
\renewcommand{\,}{\kern 0.1em}

{\newline\vspace{\abovedisplayskip}\hbox to \textwidth\bgroup\hss$\displaystyle}
{$\hss\egroup\vspace{\belowdisplayskip}}

\makeatletter
\DeclareFontFamily{OMX}{MnSymbolE}{}
\DeclareSymbolFont{MnLargeSymbols}{OMX}{MnSymbolE}{m}{n}
\SetSymbolFont{MnLargeSymbols}{bold}{OMX}{MnSymbolE}{b}{n}
\DeclareFontShape{OMX}{MnSymbolE}{m}{n}{
	<-6>  MnSymbolE5
	<6-7>  MnSymbolE6
	<7-8>  MnSymbolE7
	<8-9>  MnSymbolE8
	<9-10> MnSymbolE9
	<10-12> MnSymbolE10
	<12->   MnSymbolE12
}{}
\DeclareFontShape{OMX}{MnSymbolE}{b}{n}{
	<-6>  MnSymbolE-Bold5
	<6-7>  MnSymbolE-Bold6
	<7-8>  MnSymbolE-Bold7
	<8-9>  MnSymbolE-Bold8
	<9-10> MnSymbolE-Bold9
	<10-12> MnSymbolE-Bold10
	<12->   MnSymbolE-Bold12
}{}

\let\llangle\@undefined
\let\rrangle\@undefined
\DeclareMathDelimiter{\llangle}{\mathopen}%
{MnLargeSymbols}{'164}{MnLargeSymbols}{'164}
\DeclareMathDelimiter{\rrangle}{\mathclose}%
{MnLargeSymbols}{'171}{MnLargeSymbols}{'171}
\makeatother

\newcommand{\Pfin}[1]{\mathcal{P}_{\mathrm{fin}, 1}(#1)}
\newcommand{\Pone}[1]{\mathcal{P}_{1}(#1)}

\begin{document}
\title{Isomorphisms of reduced power monoids of semilattices}

\author{Stefano Bonzio}
\address{Department of Mathematics and Computer Science, University of Cagliari, Italy}
\email{stefano.bonzio@unica.it}

\author{Laura Cossu}
\address{Department of Mathematics and Computer Science, University of Cagliari, Italy}
\email{laura.cossu3@unica.it}



\begin{abstract}
The {\it reduced power monoid} $\Pone H$ of a monoid $H$ consists of all subsets of $H$ containing its identity, with setwise multiplication. Its submonoid of finite subsets is the {\it reduced finitary power
monoid} $\Pfin H$. For a semilattice $H$ with identity and an arbitrary monoid $K$, we prove that every isomorphism $\Pone H\to\Pone K$ or $\Pfin H\to\Pfin K$ has the form $X\mapsto\{g(x):x\in X\}$ for a
unique identity-preserving bijection $g:H\to K$. Such an isomorphism exists if and only if $K$ satisfies explicit conditions on its multiplication and $H$ and $K$ have isomorphic submonoid lattices.
Using this characterization, we identify the automorphism groups of both reduced power monoids of $H$ with the group of identity-preserving weak automorphisms of $H$ and with the automorphism group of its
submonoid lattice. Finally, we give sufficient conditions for either reduced power monoid
to determine $H$ up to isomorphism among all monoids, and show that this can fail even when $H$ and $K$ are finite semilattices with identity.
\end{abstract}

\maketitle

\thispagestyle{empty}

\section{Introduction}\label{intro}
How much of a semigroup can be recovered from the multiplication of its subsets? The classical isomorphism problem, going back to Tamura and Shafer \cite{Tam-Sha-1967}, asks when two semigroups with isomorphic
power semigroups must themselves be isomorphic. Here the power semigroup $\mathcal P(S)$ of a semigroup $S$ consists of all nonempty subsets of $S$ under setwise multiplication. Gould, Iskra, and Tsinakis
\cite[Theorem~1.3]{GIT1984} gave a positive answer when both semigroups are semilattices with identity. For arbitrary semilattices $S$ and $T$, Kobayashi \cite{Koba-1984} showed that every isomorphism
$\mathcal P(S)\to\mathcal P(T)$ maps singletons onto singletons and hence induces an isomorphism $S\to T$. For recent work on inverse semigroups and monoids, see \cite{Liu-Tri26,Dam-Loi26}.

In this paper we consider monoids, that is, semigroups with an identity element. For a monoid $H$ with identity $1_H$, the subsets of $H$ containing $1_H$ form a submonoid of $\mathcal P(H)$, called the
{\it reduced (large) power monoid} of $H$ and denoted by $\Pone H$. Its submonoid of finite subsets is the {\it reduced finitary power monoid} $\Pfin H$. In both cases, $\{1_H\}$ is the only singleton.

One motivation for studying finitary power monoids comes from additive combinatorics, where writing a finite set as a sumset amounts to factoring it in the power monoid \cite{Fa-Tr18}. For reduced finitary power monoids, arithmetic questions include atomicity \cite{An-Tr18} and minimal factorizations into irreducibles \cite{Co-Tr-25(a)}; see \cite{Tri-survey26} for a survey. 
The isomorphism problem for these monoids has been studied for numerical monoids \cite{Bie-Ger-22}, for additive submonoids of the nonnegative rationals \cite{Tri-Yan-23(a)}, and more generally
for cancellative monoids \cite{Tri-Yan2026(a),Rago26,Rago26(c)}.

One can also ask whether a given isomorphism $h$ between reduced power monoids must have the form $h(X)=\{g(x):x\in X\}$ for some identity-preserving bijection $g$ between their ground monoids. This can fail even for automorphisms. For example, the unique nonidentity automorphism of $\mathcal P_{\mathrm{fin},0}(\mathbb N_0)$ is the reversion $X\mapsto\max(X)-X$ \cite{Byr-Llo-Ped-Ste-1984,Tri-Yan2025(b)}.
It fixes every pair $\{0,n\}$ but sends $\{0,1,3\}$ to $\{0,2,3\}$. Thus any identity-preserving bijection inducing reversion would have to fix every $n$, a contradiction. Automorphisms of reduced
power monoids have also been studied for numerical monoids \cite{Bho-Tri26,Tri-Wen-2026(b)} and abelian groups
\cite{Rago26(b),Peng-Deng26}.

Our main results concern semilattices with identity, that is, commutative monoids in which every element is idempotent. We treat the large and finitary cases together, writing $\mathcal Q$ for $\mathcal P_1$ or $\mathcal P_{\mathrm{fin},1}$, with the same choice in the domain and codomain of each map. For a fixed semilattice $H$ with identity, we describe all isomorphisms $\mathcal Q(H)\to\mathcal Q(K)$ and determine which monoids $K$ can occur. Theorem~\ref{thm: reconstruction} shows that every such
isomorphism $h$ has the form
\[
h(X)=g[X]=\{g(x):x\in X\}\quad \text{ for every }X\in\mathcal Q(H)
\]
for a unique identity-preserving bijection $g:H\to K$. It also gives a necessary and sufficient condition on $g$ for this formula to define an isomorphism. The key step is to recover $h(X)$ for every $X\in\mathcal Q(H)$ from the images of the two-element sets $\{1_H,x\}$ with $x\ne 1_H$.

Although $g$ need not be a monoid homomorphism, it sends the submonoids of $H$ onto those of $K$. This leads to our characterization of the possible monoids $K$ (Theorem~\ref{thm: codomains of fixed semilattice}). We have $\mathcal Q(H)\cong\mathcal Q(K)$ if and only if the lattices of the submonoids of $H$ and $K$ are isomorphic and $K$ has the following structure. Its group of units has order at most two, and
its nonunits, if any, form a band $B$ (i.e., an idempotent semigroup). Moreover, $ub=bu=b$ for every
unit $u$ and every $b\in B$, and $\{ab,ba\}=\{a,b\}$ whenever $a,b\in B$ do not commute.
Under these conditions, the isomorphisms of the submonoid lattices are in canonical bijection with the isomorphisms $\mathcal Q(H)\to\mathcal Q(K)$ for either choice of $\mathcal Q$ (Corollary~\ref{cor: correspondence of isomorphisms}). 

For semilattices generated by their atoms, we can say more about the inducing bijection. If $H$ has at least two atoms, every isomorphism $\mathcal Q(H)\to\mathcal Q(K)$ is induced by a unique monoid isomorphism $H\to K$ (Proposition~\ref{prop: atomistic rigidity}). Thus either reduced power monoid determines such an $H$ up to isomorphism among all monoids.

Finally, we assume that $K$ is also a semilattice with identity. The automorphism groups of $\Pone H$ and $\Pfin H$ are canonically isomorphic both to the automorphism group of the submonoid lattice of $H$ and to the group of identity-preserving weak automorphisms of $H$ (Theorem~\ref{thm: automorphism groups}). When $H$ is a chain, every permutation of $H$ fixing $1_H$ induces an automorphism of $\mathcal Q(H)$. In this case, $\mathcal Q(H)\cong\mathcal Q(K)$ if and only if $K$ is a chain of the same cardinality as $H$
(Corollary~\ref{cor: chain reconstruction}). Thus finite chains are determined up to isomorphism, whereas infinite chains need not be. For semilattices $H$ and $K$ with identity, the implication $\mathcal Q(H)\cong\mathcal Q(K)\Rightarrow H\cong K$ holds whenever $|H|\le 4$ (Proposition~\ref{prop: small semilattices}), and a pair
of five-element semilattices shows that the bound is sharp.

\section{Preliminaries and elementary tools}
In this section, we collect several preliminary results, including some standard facts and a few elementary observations, that will be used throughout the paper.
Let $H$ be a multiplicatively written monoid with identity $1_H$. Submonoids of $H$ are meant to contain $1_H$. For $X\subseteq H$, we denote by $\langle X\rangle_H$ the submonoid of $H$ generated by $X$, omitting the subscript when the ambient monoid is clear. For $X=\{x_1,\ldots,x_n\}$, we also write $\langle x_1,\ldots,x_n\rangle_H$. We say that $H$ is reduced if its group of units, denoted by $H^\times$, is trivial. We refer to \cite{ClifPres} for background on semigroup theory.

In the notation of the introduction, set
\[
\Pone H=\{X\subseteq H:1_H\in X\} \;\text{ and }\; \Pfin H=\{X\subseteq H:1_H\in X,\ |X|<\infty\},
\]
both endowed with setwise multiplication.
For each $x\in H$, we write $I_x^H:=\{1_H,x\}$, omitting the superscript when $H$ is clear. Given a map $g:A\to B$ and a subset $X\subseteq A$, we write $g[X]:=\{g(x):x\in X\}$.

The next result is folklore. We prove it for completeness.
\begin{lemma}\label{lem: idempotents = submonoids}
Let $H$ be a monoid and $E\in\mathcal Q(H)$. Then $E$ is idempotent if and only if it is a submonoid of $H$. If $E$ is such a submonoid of $H$ and $X\in\mathcal Q(H)$, then 
\[XE=E \; \Longleftrightarrow \; X \subseteq E \; \Longleftrightarrow \; EX=E.\]
\end{lemma}
\begin{proof} As $1_H\in E$, the condition of idempotency $E^2=E$ is equivalent to closure under multiplication. This settles the first statement. Let $E$ be a submonoid of $H$ and $1_H\in X$. If $X\subseteq E$, then $E=\{1_H\}E\subseteq XE\subseteq E^2=E$. Conversely, if $XE=E$, then $X=X\{1_H\}\subseteq XE=E$. The right-sided argument is identical.
\end{proof}

We will make extensive use of the crucial two-to-two lemma of Tringali and Yan \cite[Corollary~3.3]{Tri-Yan2026(a)}, in the following formulation, valid for either choice of $\mathcal Q$.

\begin{lemma}[Two-to-two lemma]\label{lem: 2-to-2 lemma}
Let $H,K$ be monoids and let $h:\mathcal Q(H)\to\mathcal Q(K)$ be a monoid isomorphism.
There is a unique bijection $g:H\to K$ such that $g(1_H)=1_K$ and 
\[
h(I_x^H)=I_{g(x)}^K\;\ \text{for every }x\in H.
\]
\end{lemma}
\begin{proof}
For $\mathcal Q=\mathcal P_{\mathrm{fin},1}$, this is the original result \cite[Corollary~3.3]{Tri-Yan2026(a)}. For $\mathcal Q=\mathcal P_1$, recall from \cite[Lemma 2.2]{Tri-Wen-2026(b)} (see also \cite[Corollary 5.2]{Li-Tri25}) that the equation $AX=X^3$ has finitely many solutions $A\in \mathcal P_1(H)$ if and only if $X\in \Pfin H$. If we set $\mathcal C_H(X)=\{A\in\Pone H: AX=X^3\}$, the monoid isomorphism $h: \Pone H\to \Pone K$ induces a bijection
\[
\mathcal C_H(X)\longrightarrow\mathcal C_K(h(X)),\]
and $X\in \Pfin H$ if and only if $h(X)\in \Pfin H$. Thus, we can apply the finitary result to the restriction $h|_{\Pfin H}$. Note that this restriction property is also established in \cite[Theorem 2.3]{Tri-Wen-2026(b)}, where it is stated for cancellative monoids for the purpose of that paper.
\end{proof}

We now use the two preceding lemmas to study how isomorphisms of reduced power monoids act on submonoids and their subsets. 

\begin{lemma}\label{lem: h on submonoids}
  Let $H$ and $K$ be monoids, $h:\mathcal Q(H)\longrightarrow \mathcal{Q}(K)$ a monoid isomorphism, and $g:H\to K$ the identity-preserving bijection associated to $h$ by Lemma \ref{lem: 2-to-2 lemma}. Then, for every submonoid $E$ of $H$ belonging to $\mathcal{Q}(H)$,
  \begin{equation}\label{eq: h(E) = g[E]}
        h(E)=g[E].
\end{equation}
Moreover, if $X\in \mathcal{Q}(H)$ and $X\subseteq E$, then 
\[h(X)\subseteq g[E].\]
\end{lemma}
\begin{proof}
    Let $E\in\mathcal Q(H)$ be a submonoid of $H$. Since $h$ preserves idempotents, Lemma~\ref{lem: idempotents = submonoids} implies that  $h(E)\in \mathcal{Q}(K)$ is a submonoid of $K$. For each $x\in H$, Lemmas~\ref{lem: idempotents = submonoids} and~\ref{lem: 2-to-2 lemma} give, as in the proof of \cite[Proposition~2.7]{Tri-Wen-2026(b)},
    \[x\in E\;\ \Longleftrightarrow\;\ I_x^HE=E \;\ \Longleftrightarrow \;\ I_{g(x)}^Kh(E)=h(E)\;\ \Longleftrightarrow \;\ g(x)\in h(E).\]
    Since $g$ is surjective, these equivalences imply $h(E)=g[E]$.

    To prove the second part of the statement, assume that $X\in\mathcal Q(H)$ is a subset of $E$. Then, by Lemma \ref{lem: idempotents = submonoids}, $XE=E$ and so $h(X)h(E)=h(E)$. A further application of the same lemma gives the required inclusion.
\end{proof}

The equality $h(E)=g[E]$ can fail for an element of $\mathcal{Q}(H)$ that is not a submonoid of $H$. We will show an example in Remark~\ref{rem: general reconstruction fails}.

\medskip 

The final lemma of this section gives a criterion for an identity-preserving bijection of the underlying sets to induce an isomorphism of reduced power monoids.

\begin{lemma}\label{lem: augmentation}
Let $H$ and $K$ be monoids and $g:H\to K$ an identity-preserving bijection. The following conditions are equivalent:
\begin{enumerate}[label=(\roman*)]
\item\label{lem: augmentation1} $g_*:\Pfin H\to\Pfin K$, $X\mapsto g[X]$, is a monoid isomorphism;
\item\label{lem: augmentation2} $g_*:\Pone H\to\Pone K$, $X\mapsto g[X]$, is a monoid isomorphism;
\item\label{lem: augmentation3} for every $x,y\in H$,
\begin{equation}\label{eq: condizione h iso}\tag{C}
\{1_K,g(x),g(y),g(xy)\}
=
\{1_K,g(x),g(y),g(x)g(y)\}.
\end{equation}
\end{enumerate}
\end{lemma}
\begin{proof}
Either of the first two conditions implies \eqref{eq: condizione h iso} by multiplicativity on $I_x^HI_y^H=\{1_H,x,y,xy\}$. 

Conversely, assume condition \eqref{eq: condizione h iso} and take $X,Y\in \mathcal Q (H)$. Since both sets contain the identity,
    \[XY=\bigcup_{x\in X,\,y\in Y}I_x^HI_y^H=\bigcup_{x\in X,\,y\in Y}\{1_H,x,y,xy\}.\]
By \eqref{eq: condizione h iso},
    \[g_*(XY)=g[XY]=\bigcup_{x\in X,\,y\in Y}\{1_K,g(x),g(y),g(x)g(y)\}=\bigcup_{x\in X,\,y\in Y}I_{g(x)}^K I_{g(y)}^K=g[X]g[Y]=g_*(X)g_*(Y).\]
The map $g_*$ is also bijective with inverse $(g^{-1})_*$, and so it is an isomorphism for either choice of $\mathcal Q$.
\end{proof}

Throughout the paper, for an identity-preserving bijection $g:H\to K$, we write $g_*:\mathcal Q(H)\to\mathcal Q(K)$ for the map $X\mapsto g[X]$, for either choice of $\mathcal Q$.

\begin{remark}\label{rem: general reconstruction fails}
Not every isomorphism of reduced power monoids is of the form $g_*$ described in Lemma~\ref{lem: augmentation}. For example, let $H = \{1, a, b, c\}$ be the Klein four-group. The map $h$ interchanging $\{1, a, b\}$ and $\{1, a, c\}$ and fixing all other elements of $\mathcal Q(H)$ is an automorphism of $\mathcal Q(H)$ by \cite[Example~2]{Rago26(b)}. If $h=g_*$ for some identity-preserving bijection $g:H\to H$, then 
\[\{1,g(x)\}=h(\{1,x\})=\{1,x\} \qquad\text{ for every }x\in H.
\]
This implies that $g=\mathrm{id}_H$, contradicting that $h\ne \mathrm{id}_{\mathcal Q(H)}$. 

\medskip

Incidentally, this example also shows that equality \eqref{eq: h(E) = g[E]} in Lemma \ref{lem: h on submonoids} does not extend to non-idempotent elements of $\mathcal{Q}(H)$. For either set $X$ interchanged by $h$, we have $X^2=H$ and $h(X)\ne  X = g[X]$.
\end{remark}

\section{Isomorphisms from reduced power monoids of semilattices}\label{sec: isomorphisms from H semilattice}

\begin{center}
\textit{From now on, unless stated otherwise, $H$ denotes a semilattice with identity.}
\end{center}

In this section, we characterize the isomorphisms $\mathcal Q(H)\to\mathcal Q(K)$ for an arbitrary monoid $K$.

We order the semilattice $H$ by
\[
 x\le y \;\ \Longleftrightarrow\;\ xy=y,
\]
so that multiplication is the join and $1_H$ is the least element. We write $x\parallel y$ when $x$ and $y$ are incomparable.

\begin{remark}\label{rem: fg submonoids of a semilattice are finite}
    Every finitely generated submonoid of $H$ is a finite subsemilattice. In fact, by commutativity and idempotence, every element of a submonoid generated by $x_1,\ldots,x_n$ is a product $\prod_{i\in I}x_i$ for some $I\subseteq\{1,\ldots,n\}$, with the empty product equal to $1_H$. Thus the submonoid has at most $2^n$ elements.
\end{remark}

The next lemma is in the same spirit as \cite[Lemma~3.2]{Tri-Wen-2026(b)}: for any subset $W$ of $H$ containing $1_H$, we determine whether $x\in W$ by comparing the products of $W$ with two submonoids of $H$.
\begin{lemma}\label{lem: separating submonoids}
For $x\in H\setminus\{1_H\}$, set
\[
A_x=\{1_H\}\cup\{z\in H:z\not\le x\},
\qquad B_x=A_x\cup\{x\}.
\]
Then $A_x$ and $B_x$ are submonoids of $H$, and for every subset $W\subseteq H$ containing $1_H$,
\[
WA_x=W\cup A_x,\qquad WB_x=W\cup B_x.
\]
In particular, $x\in W$ if and only if $WA_x=WB_x$.
\end{lemma}
\begin{proof}
If $z\not\le x$ and $t\in H$, then $zt\not\le x$, as otherwise $z\le zt\le x$. This shows that $A_x$ is closed under multiplication. Moreover, if $a\in A_x\setminus\{1_H\}$, $xa\in A_x$, and $x^2=x$. Hence, $B_x$ is also closed under multiplication. Both sets contain $1_H$, so they are submonoids of $H$.

Now, let $W$ be a subset of $H$ containing $1_H$. Since $W$ and $A_x$ both contain $1_H$, we have $W\cup A_x\subseteq WA_x$. For the reverse inclusion, let $w\in W$ and $a\in A_x$. If $a=1_H$, then $wa=w$; otherwise
$a\not\le x$, so $wa\in A_x$ by the preceding observation. Thus $WA_x=W\cup A_x$.

Similarly, $W\cup B_x\subseteq WB_x$. As $B_x=A_x\cup\{x\}$, it remains to consider the products $wx$ with $w\in W$. If $w\le x$, then $wx=x$; otherwise $wx\in A_x$. This proves $WB_x=W\cup B_x=WA_x\cup\{x\}$. 

Hence $WA_x=WB_x$ if and only if $x\in W\cup A_x$, which is equivalent to $x\in W$ since $x\notin A_x$.
\end{proof}

We can now give a complete description of the isomorphisms under consideration.
\begin{theorem}\label{thm: reconstruction}
Let $H$ be a semilattice with identity and let $K$ be a monoid. A map $h:\mathcal Q(H)\to\mathcal Q(K)$ is a monoid isomorphism if and only if there exists an identity-preserving bijection $g:H\to K$ satisfying condition \eqref{eq: condizione h iso} such that
\[h(X)=g[X]=\{g(x):x\in X\}\qquad\text{for every }X\in\mathcal Q(H).\]
When these conditions hold, $g$ is unique.
\end{theorem}
\begin{proof}
Sufficiency follows from Lemma~\ref{lem: augmentation}, so we are left to prove the left to right implication.
Assume that $h:\mathcal Q(H)\to \mathcal Q(K)$ is a monoid isomorphism and let $g:H\to K$ be the unique identity-preserving bijection associated to $h$ by the two-to-two lemma \ref{lem: 2-to-2 lemma}. For $x,y\in H$, the set $E=I_x^HI_y^H=\{1_H,x,y,xy\}$ is a finite submonoid of $H$ ($H$ is commutative and idempotent). Then, Lemma~\ref{lem: h on submonoids} gives
\[
g[E]=h(E)=h(I_x^H)h(I_y^H)=I_{g(x)}^KI_{g(y)}^K.
\]
Hence $g$ satisfies condition \eqref{eq: condizione h iso}, and so $g_*:\mathcal Q(H)\to\mathcal Q(K)$ is an isomorphism by Lemma~\ref{lem: augmentation}.

We claim that $h = g_\ast$. To this end, we show that
\[\alpha:=g_*^{-1}\circ h\]
is the identity on $\mathcal Q(H)$.

By Lemma~\ref{lem: h on submonoids}, $\alpha$ fixes every submonoid $E$ of $H$ belonging to $\mathcal Q(H)$: $\alpha (E)=g^{-1}[h(E)]=g^{-1}[g[E]]=E$. 
Fix $X\in\mathcal Q(H)$ and let $E=\langle X\rangle$. Then $E\in \mathcal Q(H)$. In fact, if $\mathcal Q=\mathcal P_{\mathrm{fin},1}$, $E$ is finite by Remark~\ref{rem: fg submonoids of a semilattice are finite}.
Since $X\subseteq E$, Lemma~\ref{lem: idempotents = submonoids} gives $XE=E$. Applying $\alpha$ yields $\alpha(X)E=E$, so $\alpha(X)\subseteq E$.

For $x\in E\setminus\{1_H\}$, apply Lemma~\ref{lem: separating submonoids} to the semilattice $E$. The resulting submonoids $A_x$ and $B_x$ belong to $\mathcal Q(H)$ and are therefore fixed by $\alpha$. The same lemma, together with the injectivity of $\alpha$, yields
\[
x\in X\; \Longleftrightarrow\; XA_x=XB_x\; \Longleftrightarrow\; \alpha(X)A_x=\alpha(X)B_x\; \Longleftrightarrow\; x\in\alpha(X).
\]
Both $X$ and $\alpha(X)$ contain $1_H$ and are contained in $E$, so $\alpha(X)=X$.
This shows that $h(X)=g[X]$ for every $X\in\mathcal{Q}(H)$. Finally, any identity-preserving bijection $f:H\to K$ inducing $h$ satisfies $I_{f(x)}^K=h(I_x^H)=I_{g(x)}^K$ for every $x\in H$, so $f=g$.
\end{proof}

The next corollary is a straightforward consequence of the above theorem.

\begin{corollary}\label{cor: h preserves cardinality}
   Let $H$ be a semilattice with identity, $K$ be a monoid, and $h:\mathcal Q(H)\to\mathcal Q(K)$ be an isomorphism. For all $X,Y\in \mathcal Q(H)$,
    \[|h(X)|=|X|, \;\  \ X\subseteq Y\Longleftrightarrow h(X)\subseteq h(Y),\]
    and 
    \[h(X\cup Y)=h(X)\cup h(Y), \;\  \ h(X\cap Y)=h(X)\cap h(Y).\]
If $\mathcal Q=\mathcal P_1$, $h$ also preserves unions and intersections of arbitrary nonempty families.
If two such isomorphisms $\mathcal{Q}(H)\to \mathcal{Q}(K)$ agree on the sets $I^H_x$, they agree everywhere.
\end{corollary}

\begin{proof}
By Theorem~\ref{thm: reconstruction}, $h(X)=g[X]$ for every $X\in\mathcal Q(H)$, where $g:H\to K$ is an identity-preserving bijection. The assertions about cardinality, inclusion, unions, and intersections follow immediately. Finally, the values of $h$ on the sets $I_x^H$ determine $g$, since $h(I_x^H)=I_{g(x)}^K$, and hence determine $h$.
\end{proof}

The common formulation of Theorem \ref{thm: reconstruction} for reduced and reduced finitary power monoids allows us to give an exact relationship between the isomorphisms over $\mathcal{P}_{\rm fin,1}$ and $\mathcal P_1$.

\begin{corollary}\label{cor: unique extension}
Let $H$ be a semilattice with identity and let $K$ be a monoid. There are canonical bijections
{\small\[
\mathrm{Iso}(\Pfin H,\Pfin K)
\;\longleftrightarrow\;
\left\{
\begin{array}{@{}c@{}}
g:H\to K \text{ an identity-preserving bijection}\\
\text{satisfying \eqref{eq: condizione h iso}}
\end{array}
\right\}
\;\longleftrightarrow\;
\mathrm{Iso}(\Pone H,\Pone K).
\]}
\normalsize
On either side, the isomorphism corresponding to $g$ is given by $X\mapsto g[X]$. In particular, every isomorphism $\Pfin H\to\Pfin K$ extends uniquely to an isomorphism $\Pone H\to\Pone K$, and every isomorphism of the latter kind arises in this way.
\end{corollary}

\section{Codomains of isomorphisms from reduced power monoids of semilattices}\label{sec: codomains}
Having determined the form of the isomorphisms from reduced power monoids of semilattices, we now turn to their possible codomains. More precisely, given a semilattice $H$ with identity, we investigate the algebraic structure of the monoids $K$ for which $\mathcal Q(H)\cong\mathcal Q(K)$.

\medskip

Recall that a band is an idempotent semigroup. A band is {\it rectangular} if $aba=a$ for all $a,b$. A semigroup is {\it breakable} if $ab\in\{a,b\}$ for all $a,b$; see \cite[Definition~3.1]{Co-Tr-25(a)}. For disjoint semigroups $G,B$, the \textit{trivial ideal extension} $G\oslash B$ is the semigroup with underlying set $G\cup B$, obtained by preserving the internal multiplications and setting $ub=bu=b$ for $u\in G,b\in B$; see
\cite[Section~3, before Remark~3.3]{Co-Tr-25(a)}.

\subsection{The algebraic structure of a possible codomain}
The next proposition gives necessary conditions on $K$ for $\mathcal{Q}(K)$ to be isomorphic to $\mathcal{Q}(H)$, when $H$ is a semilattice.

\begin{proposition}\label{prop: target trivial extension}
Let $H$ be a semilattice with identity and let $K$ be a monoid such that $\mathcal Q(H)\cong\mathcal Q(K)$. Set $G=K^\times$ and $B=K\setminus G$. Then $|G|\le 2$. If $B\ne\varnothing$, then $K=G\oslash B$, where $B$ is a band such that, for all $a,b\in B$,
\begin{equation}\label{eq: band condition}\tag{B}
ab=ba\quad\text{or}\quad \{ab,ba\}=\{a,b\}.
\end{equation}
\end{proposition}
\begin{proof}
Fix an isomorphism $h:\mathcal Q(H)\to\mathcal Q(K)$, and let $g:H\to K$ be the bijection associated with $h$ by Theorem~\ref{thm: reconstruction}. For each $a\in K$, the set $I_a^K=h(I_{g^{-1}(a)}^H)$ is idempotent, since every $I_x^H$ is idempotent. Thus
\[
\{1_K,a,a^2\}=(I_a^K)^2=I_a^K=\{1_K,a\}.
\]
It follows that every unit of $K$ has square $1_K$ (since $1_K$ is the only idempotent unit), while every nonunit is idempotent.

We first show that $|G|\le2$. Suppose that $u,v\in G\setminus\{1_K\}$ are distinct, and write
\[
p=g^{-1}(u),\qquad q=g^{-1}(v).
\]
In $G$, the equality $uv=1_K$ would give $u=v^{-1}=v$, and $uv\in \{u,v\}$ would give either $u=1_K$ or $v=1_K$. Thus, $uv\notin\{1_K,u,v\}$, and condition \eqref{eq: condizione h iso} gives $g(pq)=uv$. Since $p^2=p$ and $u^2=1_K$, we have $p(pq)=pq$ and $g(p)g(pq)=u(uv)=v$. Taking $x=p$ and $y=pq$ in \eqref{eq: condizione h iso} therefore gives
\[
\{1_K,u,uv\}=\{1_K,u,uv,v\},
\]
contradicting that $v\notin\{1_K,u,uv\}$. Thus $|G|\le2$.

If $B$ is empty, the conclusion follows. Otherwise, let $b,c\in B$.
If $bc\in G$, then $b(bc)=bc$ implies $b=1_K$, a contradiction. Thus $B$ is closed under multiplication and its elements are idempotent, so $B$ is a band. 

As $H$ is commutative, so is $\mathcal Q(H)$ and then $\mathcal Q(K)$. Therefore, 
\begin{equation}\label{eq: target augmented commutativity}
\{1_K,a,b,ab\}=\{1_K,a,b,ba\}
\qquad \text{for every }a,b\in K.
\end{equation}
If $a,b\in B$ and $ab\ne ba$, this equality forces $ab,ba\in\{1_K,a,b\}$. Since both products belong to $B$, neither equals
$1_K$, and therefore
\[
\{ab,ba\}=\{a,b\}.
\]
This proves \eqref{eq: band condition}.

Finally, for $u\in G$ and $b\in B$, neither $ub$ nor $bu$ is a unit (in fact, if $ub\in G$ then also $u(ub) = b\in G$).
If $ub\ne bu$, \eqref{eq: target augmented commutativity} implies that both products must lie in $B\cap\{1_K, u,b\}=\{b\}$, which is impossible. Thus $ub=bu$. Since $u^2=1$, $b^2=b$, and $ub$ is an idempotent,
\[
ub=(ub)^2=u^2b^2=b.
\]
Therefore, $K=G\oslash B$, as required.
\end{proof}

\begin{example}\label{ex: structural conditions not sufficient}
The structural conditions on $K$ in Proposition~\ref{prop: target trivial extension} are not sufficient to guarantee that $\mathcal Q(H)\cong \mathcal Q(K)$ for a fixed semilattice $H$, even when $|H|=|K|$.

Indeed, let $H$ be a four-element chain with identity, and let $K=\{1_K,a,b,t\}$ be the four-element semilattice determined by
\[
1_K<a<t,\qquad 1_K<b<t,\qquad a\parallel b.
\]
Then
\[
K=\{1_K\}\oslash\{a,b,t\},
\]
where $\{a,b,t\}$ is a commutative band and therefore satisfies
condition~\eqref{eq: band condition}.

Since $H$ is a chain, the product of any two elements of $H$ is their maximum and hence equals one of the two factors. Thus $UV\subseteq U\cup V$ for all $U,V\in\mathcal Q(H)$, while the reverse inclusion follows from $1_H\in U\cap V$. In particular, every element of $\mathcal Q(H)$ is idempotent. In $K$, however, $ab=t$, so $X=\{1_K,a,b\}$ satisfies
\[
X^2=X\cup\{t\}\ne X.
\]
Therefore $\mathcal Q(H)$ and $\mathcal Q(K)$ are not isomorphic.
\end{example}

\begin{remark}[Hall's structural description]\label{rem: Hall structure}
A band satisfying \eqref{eq: band condition} is almost commutative in Hall's sense \cite{Hall72}:
for every $a,b$, either $ab=ba$, or $aba=a$ and $bab=b$.
By \cite[Theorem~1]{Hall72}, with our order convention, the bands
satisfying \eqref{eq: band condition} are exactly those admitting a decomposition
\[
B=\coprod_{\alpha\in Y}B_\alpha,
\]
where $Y$ is a semilattice, each $B_\alpha$ is a breakable rectangular band, and $B_{\alpha\beta}$ is a singleton whenever $\alpha\parallel\beta$. For $a\in B_\alpha,b\in B_\beta$ with $\alpha\ne\beta$,
\[
ab=ba=
\begin{cases}
b,&\alpha<\beta,\\
a,&\beta<\alpha,\\
e_{\alpha\beta},&\alpha\parallel\beta,
\end{cases}
\]
where $B_{\alpha\beta}=\{e_{\alpha\beta}\}$ in the last case. To see the specialization, note that in a rectangular band two commuting elements must coincide. In fact, if $ab=ba$, then $a=aba=ab=bab=b$. Thus, for distinct elements $a,b$ of the same rectangular component, the first alternative in \eqref{eq: band condition} is impossible, so $\{ab,ba\}=\{a,b\}$. Hence each rectangular component is breakable. Vice versa, in a breakable component, $ab,ba\in\{a,b\}$. Hence either $ab=ba$ or the two products are distinct and $\{ab,ba\}=\{a,b\}$. Elements in distinct components commute by Hall's theorem, so \eqref{eq: band condition} holds throughout $B$.
\end{remark}

\subsection{A complete criterion for codomains}
As Example~\ref{ex: structural conditions not sufficient} shows, the structural conditions on $K$ in
Proposition~\ref{prop: target trivial extension} do not suffice to ensure that $\mathcal Q(H)\cong\mathcal Q(K)$ for a fixed semilattice $H$. The remaining condition concerns the submonoid lattices of $H$ and $K$.

\medskip

For a monoid $M$, let $\mathcal L_1(M)$ denote the lattice of all submonoids of $M$, ordered by inclusion. Its least element is $\{1_M\}$. The meet of two submonoids is their intersection, and their join is the submonoid generated by their union.

The next theorem shows that, under the structural conditions on $K$, the reduced power monoids of $H$ and $K$ are isomorphic precisely when their submonoid lattices are isomorphic.

\begin{theorem}\label{thm: codomains of fixed semilattice}
Let $H$ be a semilattice with identity and let $K$ be a monoid. Then $\mathcal Q(H)\cong\mathcal Q(K)$
if and only if both of the following conditions hold:
\begin{enumerate}[label=(\roman*)]
\item\label{thm: codomains 1} $K$ satisfies the structural conditions of Proposition~\ref{prop: target trivial extension};
\item\label{thm: codomains 2} $\mathcal L_1(H)\cong\mathcal L_1(K)$.
\end{enumerate}
\end{theorem}
\begin{proof}
Suppose first that $h:\mathcal Q(H)\to\mathcal Q(K)$ is an isomorphism.
Proposition \ref{prop: target trivial extension} gives \ref{thm: codomains 1}. By Theorem~\ref{thm: reconstruction}, there is an identity-preserving bijection $g:H\to K$ inducing $h$ and satisfying condition~\eqref{eq: condizione h iso}.

Let $E\subseteq H$ be a subset containing $1_H$. For $x,y\in E$, \eqref{eq: condizione h iso} gives the following equivalences:
\[xy\in E\;\ \Longleftrightarrow\;\ \{1_K,g(x),g(y),g(xy)\}=\{1_K,g(x),g(y),g(x)g(y)\}\;\ \subseteq g[E]\;\ \Longleftrightarrow \;\ g(x)g(y)\in g[E],\]
and this applies independently of the cardinality of $E$. Thus $E$ is a submonoid of $H$ if and only if $g[E]$ is a submonoid of $K$. Consequently, $E\mapsto g[E]$ defines a lattice isomorphism $\mathcal L_1(H)\to\mathcal L_1(K)$, proving \ref{thm: codomains 2}.

Conversely, assume \ref{thm: codomains 1} and \ref{thm: codomains 2} and let $\Phi:\mathcal L_1(H)\to\mathcal L_1(K)$ be a lattice isomorphism. Every $I_a^K$ is a submonoid: for $a\in G=K^\times$, this follows from $|G|\le2$, while every nonunit is idempotent. Hence the atoms of the lattice $\mathcal L_1(K)$, that is, the submonoids minimal among those properly containing $\{1_K\}$, are precisely the
sets $I_a^K$ with $a\ne1_K$. Likewise, the atoms of $\mathcal L_1(H)$ are the sets $I_x^H$ with $x\ne1_H$. Set $g(1_H)=1_K$. For each $x\ne1_H$, define $g(x)$ to be the unique nonidentity element of $\Phi(I_x^H)$. Since $\Phi$ maps the atoms of $\mathcal L_1(H)$ bijectively onto those of $\mathcal L_1(K)$, the
resulting map $g:H\to K$ is a bijection. By construction, $\Phi(I_x^H)=I_{g(x)}^K$ for every $x\in H$.

For every submonoid $E$ of $H$,
\[x\in E \;\ \Longleftrightarrow\;\ I_x^H\subseteq E \;\ \Longleftrightarrow\;\ I_{g(x)}^K\subseteq\Phi(E).
\]
Since $g$ is surjective, this proves $\Phi(E)=g[E]$. To prove that $g$ induces an isomorphism of reduced power monoids, we need to show that it satisfies condition \eqref{eq: condizione h iso}.

Observe that, for every $a,b\in K$, 
\[
\langle a,b\rangle_K=\{1_K,a,b,ab\}.
\]
If both elements are units, this follows from $|G|\le2$; if exactly one is a unit, this follows from the trivial-extension multiplication; if both belong to $B$, this follows from \eqref{eq: band condition}. 

Since $\Phi$ preserves joins, and $\langle x,y\rangle_H=I_x^H\vee I_y^H$, 
\[
g[\langle x,y\rangle_H] =\Phi(I_x^H\vee I_y^H) = I_{g(x)}^K\vee I_{g(y)}^K =\langle g(x),g(y)\rangle_K.
\]
Now $\langle x,y\rangle_H=\{1_H,x,y,xy\}$, and we have just shown that $\langle a,b\rangle_K=\{1_K,a,b,ab\}$ for all $a,b\in K$. Taking $a=g(x)$ and $b=g(y)$ gives
\[ \{1_K,g(x),g(y),g(xy)\}=\{1_K,g(x),g(y),g(x)g(y)\}.
\]
Thus condition~\eqref{eq: condizione h iso} holds, and Lemma~\ref{lem: augmentation} shows that $g_*$ is an isomorphism for both choices of $\mathcal Q$.
\end{proof}

\begin{remark}\label{rem: optimal conditions}
Neither condition in the preceding theorem implies the other, even when $H$ is fixed and $|K|=|H|$. To see this, consider the diamond and the chain from Example~\ref{ex: structural conditions not sufficient},
this time taking the diamond as $H$ and the chain as $K$. Then $K$ satisfies \ref{thm: codomains 1}, but the chain has eight submonoids and the diamond has only seven. Hence \ref{thm: codomains 2} fails.

For an example satisfying \ref{thm: codomains 2} but not \ref{thm: codomains 1}, keep the diamond $H=\{1,a,b,t\}$, with atoms $a,b$ and greatest element $t$. Let $K=\{1,a,b,t\}$ be the monoid of transformations of $\{0,1,2\}$ given by the identity $1$, the constant maps $a$ and $t$ with values $0$ and $1$, respectively, and the map $b$ fixing $1,2$ and sending $0$ to $1$.

With the convention $fg=f\circ g$, the products of $a$ and $b$ are
\[
\begin{array}{c|cc}
    & ab & ba \\ \hline
  H & t  & t  \\
  K & a  & t
\end{array}
\]
For every other ordered pair, the product in either monoid is one of the two factors. Thus, in both monoids, $\{1,a,b\}$ is the only subset containing the identity that is not a submonoid. The two submonoid lattices therefore coincide under the given labeling, and \ref{thm: codomains 2} holds.

However, in $K$ the nonunits $a,b$ satisfy $ab\ne ba$ and
$\{ab,ba\}=\{a,t\}\ne\{a,b\}$. Hence \eqref{eq: band condition}
fails, and so does \ref{thm: codomains 1}.

In both cases, Theorem~\ref{thm: codomains of fixed semilattice} gives $\mathcal Q(H)\not\cong\mathcal Q(K)$.
\end{remark}

\begin{corollary}\label{cor: correspondence of isomorphisms}
Let $H$ be a semilattice with identity and let $K$ be a monoid such that $\mathcal Q(H)\cong\mathcal Q(K)$.
Every lattice isomorphism $\Phi:\mathcal L_1(H)\to\mathcal L_1(K)$ determines a unique identity-preserving bijection $g:H\to K$ satisfying
\[
\Phi(E)=g[E]\qquad\text{ for every }E\in\mathcal L_1(H).
\]
For either choice of $\mathcal Q$, the correspondence
\[
\Phi\longmapsto g_*
\]
defines a bijection
\[
\mathrm{Iso}_{\mathrm{lat}}(\mathcal L_1(H),\mathcal L_1(K))
\longrightarrow
\mathrm{Iso}(\mathcal Q(H),\mathcal Q(K)).
\]
\end{corollary}
\begin{proof}
In the proof of Theorem~\ref{thm: codomains of fixed semilattice} we construct $g$ from $\Phi$ and show that $g_*$ is an isomorphism for both choices of $\mathcal Q$.
Conversely, by Theorem~\ref{thm: reconstruction}, every isomorphism $\mathcal Q(H)\to\mathcal Q(K)$ is induced by a unique identity-preserving bijection $g$ satisfying condition~\eqref{eq: condizione h iso}. As in the proof of Theorem \ref{thm: codomains of fixed semilattice}, condition \eqref{eq: condizione h iso} makes $E\mapsto g[E]$ a lattice isomorphism. These constructions are inverse to each other.
\end{proof}

\subsection{Chains and semilattices generated by their atoms}
We apply Theorem~\ref{thm: codomains of fixed semilattice} to characterize the monoids $K$ for which $\mathcal Q(H)\cong\mathcal Q(K)$ when $H$ is a chain. We then prove a rigidity result when $H$ is generated by its order-theoretic atoms and has at least two of them.

\medskip

We begin with chains. Since the product of two elements of a chain is one of the two, every chain semilattice is breakable. As Tringali and Yan observed \cite[Examples~1.2(2)]{Tri-Yan-23(a)}, multiplication in the reduced finitary power monoid of a breakable monoid is given
by union. The same argument gives $XY=X\cup Y$ for arbitrary subsets $X,Y$ containing the identity.

\begin{proposition}\label{prop: chain codomains}
Let $H$ be a chain with identity and let $K$ be a monoid. Then $\mathcal Q(H)\cong\mathcal Q(K)$ if and only if $|H|=|K|$ and either $K$ is a group of order at most two, or $K=G\oslash B$, where $G=K^\times$ has order at most two and $B$ is a breakable band. Under these conditions, every identity-preserving bijection $g:H\to K$ induces an isomorphism $g_*$.
\end{proposition}
\begin{proof}
Suppose that $\mathcal Q(H)\cong\mathcal Q(K)$. Theorem~\ref{thm: reconstruction} gives $|H|=|K|$, while
Theorem~\ref{thm: codomains of fixed semilattice} gives the structural description of $K$. If $K=G\oslash B$, we only need to show that $B$ is breakable.

Since $H$ is a chain, every element of $\mathcal Q(H)$ is idempotent, and hence so is every element of $\mathcal Q(K)$. For $a,b\in B$, the set $E=\{1_K,a,b\}$ therefore satisfies $E^2=E$. As $B$ is a subsemigroup not containing $1_K$, $ab\in E\cap B=\{a,b\}$. Thus $B$ is breakable.

Conversely, suppose that $|H|=|K|$ and that $K$ has one of the two stated forms. In either case, $ab\in\{1_K,a,b\}$ for all $a,b\in K$. Thus, for $U,V\in\mathcal Q(K)$, we have $UV\subseteq U\cup V$, and the reverse inclusion holds since both sets contain $1_K$. Hence multiplication in $\mathcal Q(K)$ is union, as it is in $\mathcal Q(H)$.

Since $|H|=|K|$, we may choose a bijection $g:H\to K$ with $g(1_H)=1_K$. For every such $g$ and all $X,Y\in\mathcal Q(H)$,
\[
g[XY]=g[X\cup Y]=g[X]\cup g[Y]=g[X]g[Y].
\]
Thus every identity-preserving bijection $g:H\to K$ induces an isomorphism $g_*$.
\end{proof}

The following example shows that the monoid $K$ in Proposition~\ref{prop: chain codomains} need not be commutative.
\begin{example}\label{ex: left zero}
Let $H=\{1<a<b\}$ be a three-element chain and let $K=\{1,u,v\}$ be the breakable monoid defined by: $u^2=u$, $v^2=v$, $uv=u$, and $vu=v$ (cf. \cite[Example 1.2(2)]{Tri-Yan-23(a)}). It is routine to check that the bijection $1\mapsto1$, $a\mapsto u$, $b\mapsto v$ induces an isomorphism $\mathcal Q(H)\to Q(K)$. Nevertheless, $K$ is not commutative.
\end{example}

We next prove a rigidity result for semilattices generated by their atoms. Here, an \textit{atom} is an element $a\in H$ covering the identity, i.e., $1_H<a$ and there is no $x\in H$ such that $1_H<x<a$.

\begin{proposition}\label{prop: atomistic rigidity}
Let $H$ be a semilattice with identity generated by its atoms, and suppose that $H$ has at least two atoms. Let $K$ be a monoid. Every isomorphism $\mathcal Q(H)\longrightarrow\mathcal Q(K)$ is induced by a unique monoid isomorphism $g:H\to K$.
In particular,
\[\mathcal Q(H)\cong\mathcal Q(K) \;\ \Longleftrightarrow \;\ H\cong K.\]
\end{proposition}

\begin{proof}
Let $h: \mathcal{Q}(H)\to \mathcal{Q}(K)$ be a monoid isomorphism. By Theorem~\ref{thm: reconstruction}, we have $h=g_*$ for a unique identity-preserving bijection $g:H\to K$ satisfying condition~\eqref{eq: condizione h iso}. Let $A$ be the set of atoms of $H$, and write $x'=g(x)$.

For distinct $a,b\in A$, the elements $a,b$ are incomparable. Thus $ab\notin\{1_H,a,b\}$, and
condition~\eqref{eq: condizione h iso}, applied to $(a,b)$ and $(b,a)$,
gives
\[
a'b'=b'a'=(ab)'.
\]
Moreover, the same condition applied to $(a,a)$ gives
$(a')^2\in\{1_K,a'\}$.
Suppose that $(a')^2=1_K$. Choose $b\in A\setminus\{a\}$, and set $c=ab$. Then
\[
a'c'=a'(a'b')=(a')^2b'=b'.
\]
Since $ac=c$, condition~\eqref{eq: condizione h iso} also gives
\[
a'c'\in\{1_K,a',c'\},
\]
a contradiction. Therefore the elements of $g[A]$ are commuting idempotents of $K$.

Let $x\ne 1_H$. Since $H$ is generated by its atoms, write $x=a_1\cdots a_n$ with each $a_i$ an atom and $n$ minimal. For $1\le i<n$, the elements $a_1\cdots a_i$ and $a_{i+1}$ are incomparable; otherwise their product would equal one of them, yielding a shorter factorization of $x$.

For incomparable $u,v\in H$, condition~\eqref{eq: condizione h iso} gives $g(u)g(v)=g(uv)=g(v)g(u)$, since $uv\notin\{1_H,u,v\}$. Applying this to the successive partial products yields
\[
g(x)=g(a_1)\cdots g(a_n).
\]
Since $g$ is surjective, $g[A]$ generates $K$ as a monoid, and $K$ is therefore a semilattice with identity.

Now that $K$ is a semilattice, condition~\eqref{eq: condizione h iso} shows that $g$ preserves comparability and incomparability in both directions. In particular, $g[A]$ is an antichain.
Moreover, each $g(a)\in g[A]$ is an atom of $K$. In fact, if $1_K<z\le g(a)$, express $z$ as a nonempty product of elements of $g[A]$. Every factor lies below $g(a)$, so the antichain property forces every factor to equal $g(a)$. Thus $z=g(a)$.

It remains to show that $g$ preserves multiplication. We start proving that $g(ax)=g(a)g(x)$ for every $a\in A$ and $x\in H$. This is immediate if $x=1_H$ and follows from \eqref{eq: condizione h iso} if $a\parallel x$. Suppose, then, that $x\ne 1_H$ is comparable with $a$. Since $a$ is an atom, $a\le x$, so $ax=x$. By \eqref{eq: condizione h iso}, $g(a)$ and $g(x)$ are comparable in $K$. As $g(a)$ is an atom and $g(x)\ne 1_K$, we have $g(a)\le g(x)$. Hence
\[
g(ax)=g(x)=g(a)g(x).
\]
Since the atoms generate $H$, induction on the length of a factorization of $x$ into atoms now gives $g(xy)=g(x)g(y)$ for all $x,y\in H$. Thus $g$ is a monoid isomorphism, and its uniqueness follows from Theorem~\ref{thm: reconstruction}.

The reverse implication is straightforward: it suffices to apply any monoid isomorphism $H\to K$ elementwise.
\end{proof}

The assumption on the number of atoms cannot be omitted: the two-element semilattice and the cyclic group of order two have isomorphic reduced power monoids; see \cite[Examples~1.2(1)]{Tri-Yan-23(a)}.

\section{Automorphism groups and the isomorphism problem}

The preceding sections describe the isomorphisms $\mathcal Q(H)\to\mathcal Q(K)$ when $H$ is a semilattice with identity,  and characterize their possible ground monoids $K$. We now assume that $K$ is also a semilattice. In this setting, condition~\eqref{eq: condizione h iso} has a classical order-theoretic interpretation in terms of weak isomorphisms. 

\medskip

We follow the terminology of \v{S}evrin \cite{Sev1965} mentioned in \cite[p.~150]{Jones1978}, using the dual order.

\begin{definition}\label{def: weak isomorphism}
Let $H,K$ be semilattices. A bijection $g:H\to K$ is a \emph{weak isomorphism} if, for every $x,y\in H$,
\begin{enumerate}[label=(\roman*)]
\item $x\parallel y$ if and only if $g(x)\parallel g(y)$, and
\item $x\parallel y$ implies $g(xy)=g(x)g(y).$
\end{enumerate}
When $H,K$ have identities, preservation of the identity is an additional condition.
\end{definition}

\begin{corollary}\label{cor: semilattice targets}
Let $H,K$ be semilattices with identity. The map $g\mapsto g_*$ is a bijection between identity-preserving
weak isomorphisms $H\to K$ and monoid isomorphisms $\mathcal Q(H)\to\mathcal Q(K)$.
\end{corollary}
\begin{proof}
By Theorem~\ref{thm: reconstruction}, it suffices to show that an identity-preserving bijection $g:H\to K$ satisfies \eqref{eq: condizione h iso} if and only if it is a weak isomorphism.

If $x,y$ are comparable, then $xy\in\{x,y\}$, so \eqref{eq: condizione h iso} holds for $x,y$ precisely when
$g(x)$ and $g(y)$ are comparable. If $x\parallel y$, then $xy\notin\{1_H,x,y\}$, so \eqref{eq: condizione h iso} holds precisely when $g(x)g(y)=g(xy)$. In this case the product differs from both $g(x)$ and $g(y)$, so $g(x)\parallel g(y)$. Thus \eqref{eq: condizione h iso} is equivalent to preserving comparability in both directions and preserving products of incomparable elements, as required.
\end{proof}

\begin{remark}\label{rem: Jones correspondence}
Let $\Sigma'(S)$ denote the lattice of all subsemilattices of a semilattice $S$, including the empty set. 
By \cite[Lemma~1.1]{Jones1978}, every lattice isomorphism $\Psi:\Sigma'(S)\to\Sigma'(T)$ is induced by a unique weak isomorphism $f:S\to T$, in the sense that $\Psi(A)=f[A]$ for every $A\in\Sigma'(S)$. As Jones notes immediately after the proof, the converse is due to \v{S}evrin: every weak isomorphism induces
a lattice isomorphism by the same formula.

Now let $H$ and $K$ be semilattices with identity. The map $E\mapsto E\setminus\{1_H\}$ identifies $\mathcal L_1(H)$ with $\Sigma'(H\setminus\{1_H\})$, and similarly for $K$. Under these
identifications, a lattice isomorphism $\Phi:\mathcal L_1(H)\to\mathcal L_1(K)$ determines a unique weak
isomorphism $f:H\setminus\{1_H\}\to K\setminus\{1_K\}$. Extending $f$ by sending $1_H$ to $1_K$ gives the unique identity-preserving weak isomorphism $g:H\to K$ such that $\Phi(E)=g[E]$ for every
$E\in\mathcal L_1(H)$. In our setting, Corollaries~\ref{cor: semilattice targets} and~\ref{cor: correspondence of isomorphisms} yield this correspondence directly.
\end{remark}

\subsection{Automorphism groups}
Let $\mathrm{WAut}_1(H)$ denote the group of weak automorphisms of $H$ fixing $1_H$, under composition. The next result specializes Remark \ref{rem: Jones correspondence} to automorphisms of reduced power monoids of semilattices with identity.

\begin{theorem}\label{thm: automorphism groups}
Let $H$ be a semilattice with identity. The map $g\mapsto g_*$ defines a group isomorphism
\[
\mathrm{WAut}_1(H)\longrightarrow
\mathrm{Aut}(\mathcal Q(H)).
\]
Consequently,
\[
\mathrm{Aut}(\Pfin H)\cong
\mathrm{Aut}(\Pone H)\cong
\mathrm{WAut}_1(H)\cong
\mathrm{Aut}_{\rm lat}(\mathcal L_1(H)).
\]
The canonical map $\mathrm{Aut}(H)\to\mathrm{Aut}(\mathcal Q(H)), \; g\mapsto g_*$ is surjective if and only if every identity-preserving weak automorphism of $H$ is a monoid automorphism.
\end{theorem}
\begin{proof}
By Corollary~\ref{cor: semilattice targets}, the map $g\mapsto g_*$ is a bijection from $\mathrm{WAut}_1(H)$ onto $\mathrm{Aut}(\mathcal Q(H))$. By Corollary~\ref{cor: correspondence of isomorphisms}, the map
$g_*\mapsto\Phi_g$, where $\Phi_g(E)=g[E]$ for every $E\in\mathcal L_1(H)$, is a bijection from $\mathrm{Aut}(\mathcal Q(H))$ onto $\mathrm{Aut}_{\rm lat}(\mathcal L_1(H))$.
Since $(g\circ f)[X]=g[f[X]]$ for every $X\subseteq H$, we have
\[
(g\circ f)_*=g_*\circ f_*
\;\; \text{ and }\;\;
\Phi_{g\circ f}=\Phi_g\circ\Phi_f.
\]
Thus both bijections are group isomorphisms.

Finally, the canonical homomorphism is the composite
\[\mathrm{Aut}(H)
\hookrightarrow \mathrm{WAut}_1(H)\xrightarrow{\;\cong\;}\mathrm{Aut}(\mathcal Q(H)).
\]
It is therefore surjective if and only if $\mathrm{Aut}(H)=\mathrm{WAut}_1(H)$, which is precisely the stated condition.
\end{proof}

\medskip

We now specialize the above theorem to chains and to semilattices generated by their atoms.

\medskip

If $H$ is a chain, every bijection $H\to H$ fixing $1_H$ is a weak automorphism, since there are no incomparable pairs. The next result follows from this observation or, equivalently, from the fact that multiplication in $\mathcal Q(H)$ is the union.
\begin{corollary}\label{cor: chain automorphisms}
If $H$ is a chain with identity, then
\[
\mathrm{Aut}(\mathcal Q(H))
\cong\mathrm{Sym}(H\setminus\{1_H\}).
\]
For a finite chain of cardinality $n$, this group is $S_{n-1}$.
\end{corollary}
    
\begin{example}\label{ex: chain nonmultiplicative}
Let $H=\{1<a<b\}$. The map $g:H\to H$ fixing $1$ and exchanging $a$ and $b$ is a weak automorphism and hence induces an automorphism
$g_*$ of $\mathcal Q(H)$. However,
\[
g(ab)=g(b)=a\neq b=g(a)g(b),
\]
so $g$ is not a monoid homomorphism. Thus the canonical inclusion $\mathrm{Aut}(H)\longrightarrow\mathrm{Aut}(\mathcal Q(H))$ need not be surjective.
\end{example}

Using Proposition \ref{prop: atomistic rigidity}, we are able to give a clean characterization of automorphisms of $\mathcal{Q}(H)$, when $H$ is generated by its (order-theoretic) atoms. 

\begin{corollary}\label{cor: atomistic automorphisms}
Let $H$ be a semilattice with identity generated by its atoms. The canonical homomorphism
\[\mathrm{Aut}(H)\to\mathrm{Aut}(\mathcal Q(H)), \; g\mapsto g_*
\]
is a group isomorphism.
\end{corollary}
\begin{proof}
If $H$ has at least two atoms, the assertion follows from Proposition \ref{prop: atomistic rigidity} with $K = H$. Otherwise, since its atoms generate $H$, we have $|H|\le 2$, and both automorphism groups are trivial.
\end{proof}

\subsection{The isomorphism problem}\label{subsec: isomorphism problem}

In this final subsection, we study the isomorphism problem for reduced power monoids of semilattices. More precisely, we seek conditions on a semilattice $H$ with identity under which $\mathcal Q(H)\cong\mathcal Q(K)$ implies $H\cong K$ for every semilattice $K$ with identity. We also give counterexamples
showing that this implication fails in general.

\medskip

We begin with chains. Gould, Iskra, and Tsinakis \cite[Theorem~1.4]{GIT1984} showed that two chains with isomorphic power semigroups are order-isomorphic. For reduced power monoids, the following corollary shows that, among chains with identity, only the cardinality is determined.

\begin{corollary}\label{cor: chain reconstruction}
Let $H,K$ be semilattices with identity and suppose that $H$ is a chain. Then
\[\mathcal Q(H)\cong\mathcal Q(K) \;\ \Longleftrightarrow\;\
K\text{ is a chain and }|H|=|K|.
\]
If $H$ is finite, these conditions imply $H\cong K$.
\end{corollary}
\begin{proof}
By Corollary~\ref{cor: semilattice targets}, an isomorphism $\mathcal Q(H)\to\mathcal Q(K)$ induces an identity-preserving weak isomorphism $H\to K$. Thus $|H|=|K|$, and since every pair of elements
of $H$ is comparable, the same holds in $K$. Conversely, if $K$ is a chain of the same cardinality as $H$, we may choose an identity-preserving bijection $H\to K$. Every such bijection is a weak isomorphism, so the
same corollary applies. The final assertion follows since finite chains of the same cardinality are order-isomorphic.
\end{proof}

The next example shows two infinite chains with the same cardinality that are not isomorphic.
\begin{example}\label{ex: infinite chains}
Let $H=\mathbb N_0$ and $K=\mathbb Q_{\ge0}$, each with multiplication given by maximum and identity $0$. They are countably infinite chains, so $\mathcal Q(H)\cong\mathcal Q(K)$ by the above result. However, they are not isomorphic:
$H$ has a least element strictly above its identity, while $K$ does not.
\end{example}

Proposition~\ref{prop: atomistic rigidity} gives a positive answer to the isomorphism problem for semilattices generated by their atoms, provided they have at least two atoms. When $K$ is also a semilattice, the hypothesis on the number of atoms can be dropped.

\begin{corollary}
    Let $H$ and $K$ be semilattices with identity and assume that $H$ is generated by its atoms. Then
\[\mathcal Q(H)\cong\mathcal Q(K) \;\ \Longleftrightarrow\;\
H\cong K.
\]
\end{corollary}
\begin{proof}
If $H$ has at least two atoms, apply Proposition \ref{prop: atomistic rigidity}. Otherwise $|H|\le 2$. Theorem \ref{thm: reconstruction} gives $|K| = |H|$, and the unique identity-preserving bijection between these semilattices is a monoid isomorphism.
\end{proof}

We conclude by showing that the isomorphism problem has a positive answer for semilattices with identity of cardinality at most four.
\begin{proposition}\label{prop: small semilattices}
Let $H$ and $K$ be semilattices with identity such that $\mathcal Q(H) \cong \mathcal Q(K)$. If $|H| \le 4$, then $H \cong K$.
\end{proposition}
\begin{proof}
By Theorem~\ref{thm: reconstruction}, $|K|=|H|$. If $H$ is a chain, the conclusion follows from Corollary~\ref{cor: chain reconstruction}. Otherwise, choose incomparable elements $a,b\in H$. Then
$1_H,a,b,ab$ are distinct and hence exhaust $H$, since $|H|\le 4$. Thus $H$ is the four-element diamond, generated by its two atoms $a,b$, and Proposition~\ref{prop: atomistic rigidity} applies.
\end{proof}

The bound on the size of $H$ given in Proposition~\ref{prop: small semilattices} is sharp. We show in the next example two non-isomorphic semilattices $H$ and $K$ with $5$ elements with $\Pfin H=\Pone H \cong \Pone K=\Pfin K$.

\begin{example}\label{ex: counterexample finite}
Let $H=\{1,a,b,c,t\}$ and $K=\{1,u,x,y,z\}$ be the semilattices whose orders are described by
\[
1<a,b<c<t,\;\  a\parallel b;
\qquad
1<u<x,y<z,\;\  x\parallel y.
\]
Multiplication is join in both cases. Their Hasse diagrams, with the
identity at the bottom, are:

\medskip

\begin{center}
\begin{tikzpicture}[scale=0.82,every node/.style={circle,draw,inner sep=2pt}]
\node (h1) at (0,0) {$1$};
\node (ha) at (-0.7,1) {$a$};
\node (hb) at (0.7,1) {$b$};
\node (hc) at (0,2) {$c$};
\node (ht) at (0,3) {$t$};
\draw (h1)--(ha)--(hc)--(ht);
\draw (h1)--(hb)--(hc);
\node[draw=none] at (0,-0.6) {$H$};
\node (k1) at (5,0) {$1$};
\node (ku) at (5,1) {$u$};
\node (kx) at (4.3,2) {$x$};
\node (ky) at (5.7,2) {$y$};
\node (kz) at (5,3) {$z$};
\draw (k1)--(ku)--(kx)--(kz);
\draw (ku)--(ky)--(kz);
\node[draw=none] at (5,-0.6) {$K$};
\end{tikzpicture}
\end{center}
The bijection defined by
\[
g(1)=1,\;\  g(a)=x,\;\  g(b)=y,\;\  g(c)=z,\;\  g(t)=u
\]
preserves the unique incomparable pair and its product: $g(ab)=g(c)=z=xy$. It is therefore a weak isomorphism preserving the identity, and induces $\mathcal Q(H)\cong\mathcal Q(K)$.
However, $H$ has two atoms in its ground order, whereas $K$ has one; hence $H\not\cong K$.
\end{example}

\section*{AI declaration}

The authors guided successive revisions of an initial draft, using
ChatGPT (OpenAI; GPT-6, GPT-6.1, and GPT-6.2, including Sol and Astra
variants) and Claude Opus 5.5 (Anthropic) to explore examples and
generalizations, compare arguments, locate references, and refine the
exposition. They assessed the suggestions and reworked those they chose to incorporate.

\section*{Acknowledgments}

The authors are grateful to Salvatore Tringali for fruitful discussions during his visit to the University of Cagliari in June 2026, which inspired this paper.

The authors acknowledge financial support from Fondazione di Sardegna through the project LoMEA and from INdAM-GNSAGA (Gruppo Nazionale per le Strutture Algebriche, Geometriche e loro Applicazioni). S.~Bonzio also acknowledges support from the European Commission through the project SemPER (Grant Agreement ID: 101299559).

\end{document}